\let\ifiso\iffalse
\let\ifdump\iffalse
\let\ifeight\iffalse
\let\ifmac\iffalse

\catcode`\{=1 
\catcode`\}=2 
\catcode`\$=3 
\catcode`\&=4 
\catcode`\#=6 
\catcode`\^=7 \catcode`\^^K=7 
\catcode`\_=8 \catcode`\^^A=8 
\catcode`\^^I=10 
\chardef\active=13 \catcode`\~=\active 
\catcode`\^^L=\active \outer\def^^L{\par} 

\catcode`\@=11
\let\e@\expandafter
\def\s@tt#1#2{\e@\e@\e@\let\e@#1#2}
\def\letifundefined#1#2{\def\aaa{#1}\def\bbb{#2}%
\ifx#1\undefined\s@tt\aaa\bbb\fi}
\letifundefined\ifdump\iftrue
\letifundefined\ifmac}|
\read\@psffile to \h@p
\def\n@xt{|
\read\@psffile to \h@p
\e@\BBt@st\h@p :
\n@xt}|
\n@xt
}|
\closein\@psffile
\e@\h@ight\t@p
\e@\advance\e@\h@ight\e@-\b@t
\e@\w@dth\r@ght
\e@\advance\e@\w@dth\e@-\l@ft}
}

\def\illt@mpl#1{\def\t@mp{illustration #1}}
\def\illbox(#1){\p@llbox{#1}\illg@tdim\illt@mpl}

\def\illustration(#1){\vskip 6pt plus 2pt minus 2pt\centerline{\illbox(#1)}%
\vskip 6pt plus 2pt minus 2pt\noindent}

\let\psprol\relax
\let\psend\relax
\def\pst@mpl#1{\e@\w@dth\l@ft \e@\h@ight\b@t
\sc@le\w@dth\scalep@rt \sc@le\h@ight\scalep@rt
\divide\w@dth65536 \divide\h@ight65536
{\count255=\scalep@rt \divide\count255 10
\count254\w@dth \count253\h@ight
\xdef\t@mp{psfile=\namep@rt hoffset=-\number\count254 \space\space
voffset=-\number\count253\space\space vscale=\number\count255\space\space
hscale=\number\count255}%
}}

\def\psbox(#1){\psprol\p@llbox{#1}\illg@tdim\pst@mpl\psend}

\def\psillustration(#1){\vskip 6pt plus 2pt minus 2pt\centerline{\psbox(#1)}%
\vskip 6pt plus 2pt minus 2pt\noindent}

\ifmac\else
\let\illbox\psbox
\let\illustration\psillustration
\fi

\langdef\yrs\in english={Yours\iffriendly\else sincerely\fi},%
\in swedish={\iffriendly H\else V{\"a}nliga h\fi {\"a}lsningar},%
\in french={Bien {\'a} \iffriendly toi\else vous\fi}.

\def\Torsten{\vskip1cm\leftline{\hskip10cm \yrs}\nobreak \vskip.75cm\nobreak
\leftline{\hskip10cm Torsten \iffriendly\else Ekedahl\fi}}

\def\letter{\footline={\hss\ifnum\pageno=1\else\tenrm\folio\fi\hss}%
\headline={\ifnum\pageno>1T.~Ekedahl\hss\today\else\hfil\fi}%
\rightline{Stockholm \today\qquad}\vskip3cm}

\langdef\dear\in english={Dear},\in swedish={},\in french={Cher}.
\def\beginletter#1\par{\letter \dear #1\unskip,\par\qquad}

\def\endletter#1\par{\if\relax#1\relax
\else\friendlytrue\fi\nobreak\Torsten\vfill\eject\@nd}

\catcode`\@=12

\ifdump
\let\dmp\dump
\else
\let\dmp\relax
\fi
\dmp

\final
\def\Tr{\symb{Tr}}
\def\Irr{\symb{Irr}}

\mdef\coh{H^*(X,r)}
\mdef{\hodg#1.#2.}{H^{#1}(X,{\gOm}^{#2}_{X/k})}

\candef W
\head varieties of CM-type

\begin{document}

\begin{start} Varieties of CM-type \by Torsten Ekedahl
\end{start}

\begin{introduction}
We will introduce the notion of a variety (or more generally a motive) of
CM-type which generalises the well known notion of abelian variety of
CM-type. Just as in that particular case it will turn out that the
cohomology of the variety is determined by purely combinatorial data; the
type of the variety. As applications we will show that the \l-adic
representations are given by algebraic Hecke characters whose algebraic
parts are determined by the type and give a method for computing the
discriminant of the N{\'e}ron-Severi group of super-singular Fermat surfaces.
\end{introduction}

\begin{section}Preliminaries.

To begin with let us recall the following facts from category theory. If
\Cal A is an additive category all of whose idempotents have kernels and
$R$ is a ring, then for a finitely right resp.~left projective $R$-module
$P$ resp.~$Q$ and a left $R$-object $M$ in \Cal A we can define objects
$P\bigotimess _RM$ resp.~$Hom_R(Q,M)$ of \Cal A characterised by
$$
\tagno{Hom_{\Cal A}(P\bigotimess _RM,N)&=Hom_R(P,Hom_{\Cal A}(M,N))\cr
\noalign{\leftline{resp. }}
Hom_{\cal A}(N,Hom_R(Q,M))&=Hom_R(Q,Hom_{\cal A}(N,M)).\cr}
$$

We always have a natural $R$-morphism $ev\co P\to Hom_{End_R(P)-\Cal
A}(Hom_R(P,M),M)$, the evaluation map, obtained by interpreting an element
$p\in P$ as an $R$-morphism $R\to P$ and using $M=Hom_R(R,M)$.  If
R=$\bigoplus P_i^{n_i}$ and $Hom_R(P_ i,P_j)=0$ for $i\ne j$ then for any
$R$-object $M$ in \Cal A we have
$$
M=\bigoplus  P_i{\textstyle\bigotimess }_{S_i}Hom_R(P_i,M),\tag 1.1
$$
where $S_i:=End_R(P_i)$ and the map is defined using the evaluation maps.
To see this we first note that $P_i=Hom_R(R,P_i)\cong(S_i)^{n_i}$ so that
$P_i$ is $S_i$-projective and then the desired equivalence follows by
decomposing the two $R$-factors of $M=R\bigotimess _RHom_R(R,M)$.

The following setup will be with us during the rest of the paper: We let
\k\ be a perfect field of characteristic $p\ge 0$, $X$ a proper, smooth
variety over \k\ (alternatively a motive) and $S$ a set of
\k-correspondences of $X$. Furthermore, \coh, $r$ prime, will denote the
\l-adic cohomology of $X_{\k}$, where \bk\ is a fixed algebraic closure of
\k, when $r\ne p$ and the crystalline cohomology of $X/\k$ when
$r=p$. Recall that when $r\ne p$ \coh\ is a graded $\Z_r$-algebra, finitely
generated as $\Z_r$-module, having a continuous action of $Gal(\bk/\k)$ and
that when $r=p$, \coh\ is a graded $\W(\k)$-algebra, finitely generated as
$\W(\k)$-module having a {\gsi}-linear endomorphism $F$. Here $\W(\k)$ is the
ring of Witt vectors of \k\ and {\gsi} sends a Witt vector $(x_i)$ to
$(x^p_i)$. We let $L_r$ be an algebraically closed field containing $\Z_r$
resp.~$\W(\k)$. Furthermore we will denote by $K$ the
fraction field of $\W(\k)$.

Finally, if $p>0$ we will have need of the following technical condition.
There is a scheme $T$ of finite type over $\F_p$, a smooth and
proper morphism $\Cal X\to T$ and a cartesian diagram
$$
\diagram{ X&\mapright{}&\cal {X}\cr
\mapdown{}&&\mapdown{}\cr
\Sp\;k&\mapright{}&T\cr}
$$
such that for every closed point $t\in T$ the eigenvalues of the Frobenius
with respect to $\k(t)$ on $H^i(X_t,p)$ are algebraic integers all of
whose archimedean absolute values are $|k(t)|^{i/2}$.

This condition is fulfilled when $X$ (possibly over \bk) is the image of a
smooth and projective variety (\[K-M]) and that this is always the case has
recently been verified by J.~de Jong (\[Jo]).

For a field $L$ and a set $R$ let $K(R,L)$ be the Grothendieck group of the
category of finite dimensional representations (i.e.~maps of $R$ into the
set of endomorphisms) of $R$. Then $K(R,L)$ is a functor in $R$ and $L$;
contravariant in $R$ and covariant in $L$. If $M(R)$ is the free monoid
generated by $R$ then the trace map gives an additive map $\Tr: K(R,L)\to
L^{M(R)}$.
\begin{lemma}1.2.
Let $L$ be algebraically closed of characteristic 0.

i) $\Tr: K(R,L)\to L^{M(R)}$ is injective.

ii) Let $L'\subseteq  L$ be a subfield of $L$ and $N$ a semi-simple
$L$-representation of $R$ s.t.~for all $r\in M(R)$ $\Tr_N(r)\in  L'$.
If $L'\langle R\rangle$ is the free associative $L'$-algebra
on $R$ then
$I :=\ker (L'\langle R\rangle\to End_L(N))$ depends only on the function
$\Tr_N: M(R)\to L'$ and $L'\langle R\rangle/I\bigotimess
_{L'}L\to  End_L(N)$ is
an injection. In particular, if $L'$ is algebraically closed
$N$ is isomorphic to the
scalar extension of some $L'$-representation of $R$.

iii) If $L'$ is an algebraically closed field and $L'\to L$
a field homomorphism,
then the following diagram
$$
\diagram{
K(R,L')&\mapright{}&K(R,L)\cr
\mapdown{}&&\mapdown{}\cr
L'{}^{M(R)}&\mapright{}&L^{M(R)}\cr}
$$
is cartesian.
\pro Let us begin with ii). Note first that as
$\im(L\langle R\rangle\to End_L(N))$ is semi-simple and that for a finite
dimensional semi-simple L-algebra $M$ and a faithful finite dimensional
$L$-representation $V$, the linear form
$$
\eqalign{ M \times\;M&\to L\cr
                 (m,m')\mapsto  Tr_V(mm')\cr}
$$
is non-degenerate. Hence $t\in L\langle R\rangle$ acts as zero on $N$ iff
$\Tr_N(rt)=0$
for all $r\in M(R)$. If $t=\sum _{r\in M(R)}{\gla}_rr$ then this is
a set of linear conditions on the ${\gla}_r$ with coefficients in $L'$
depending only on $\Tr_N: M(R)\to L'$. Furthermore,
$L'\langle
R\rangle/I\bigotimess _{L'}L\to End_L(N)$ is injective iff whenever
there is an $L$-linear relation in $End_L(N)$ between elements in
$L'\langle R\rangle$ there is also an $L'$-linear relation. This also follows
from the fact that the above conditions have $L'$-coefficients. Now iii)
follows immediately from i) and ii) whereas i) is well known
(cf.~\[C-R:Thm. 30.12]).
\end{lemma}

  If $L$ is a field of characteristic zero and $L'$ an
algebraic closure of $L$, we
put $\ovl K(R,L):=K(R,L')\bigcap  L^{M(R)}$. This is clearly independent of the
choice of $L'$ and we have $K(R,L)\subseteq \ovl K(R,L)$.
If $N$ is a representative over $L'$ of an $n\in \ovl K(R,L)$ then we can
construct the $L\langle R\rangle/I$ of lemma \ref{1.1}, which
depends only on $n$. It is a semi-simple $L$-algebra as its
scalar extension to $L'$ is, and we will denote it $A_L(n)$. In case $L=\Q$
then we put $A(n):=\im(\Z\langle R\rangle\to A_\Q(n)$, where $\Z\langle
R\rangle$ denotes the free associative algebra on $R$.
If $M$ is an over-field of $L$ then we say that $n$ is
realisable over $M$ if the
induced element in $\ovl K(R,M)$ belongs to $K(R,M)$. This is equivalent to $N$
being realisable by an $A_L(n)\bigotimess _LM$-representation. For any
$n\in K(R,L')$ we let $\Irr(n)\subset K(R,L')$ be the set of irreducible
constituents of n. If now $n\in\ovl K(R,L)$ then $\Irr(n)$
is a finite set stable
under the action of $Gal(\bar\Q/\Q)$. Under the correspondence of Galois
theory, $\Irr(n)$ then corresponds the {\'e}tale $L$-algebra $Z(A_L(n))$.
If we return to the situation at hand we have elements
$[H^i(X,r)\bigotimess _{\Z_r}L_r]$
(resp.~$[H^i(X,p)\bigotimess _{\W(k)}L_p]$) in $K(S,L_r)$. Let
$\bar \Q$ be an algebraic closure of \Q.
\begin{lemma}1.3.  Let $K_r$ denote $\Q_r$ when $r\ne p$ and
$K_r$ when it isn't. There
exists a unique element $[H^i(X)]\in \ovl K(S,\Q)$ whose image in
$K(S,K_r)$  coincides with $[H^i(X,r)]$.
Furthermore, $A([H_i(X)])$ is finitely generated as \Z-module and
$A\bigotimess K_r$  equals $A_r/rad(A_r)$ where
$$
A_r:=Im(\Q_r\langle S\rangle\to
End_{\Q_r}(H^i(X,r))\;\;\;\;\;\;(resp.\;\ldots).
$$
\pro I first claim that, for every $s\in M(S)$, $Tr(s,H_i(X,r))$ is a rational
number independent of $r$. By standard specialisation arguments we reduce to
\k\ being a finite field where it is \[K-M:Thm 2]
(supplemented by \[Gr] for the
definition of the cycle map in crystalline cohomology).
Note that if $p=0$, using Chow's lemma and resolution of singularities
we can get a reduction for which our technical condition is fulfilled. In this
case, a transcendental argument can also be used. This already,
using \(1.2:ii),
gives the existence of  $[H_i(X)]$ and that $A\bigotimess L_r
= A_r/rad(A_r)$. Hence  $A\bigotimess \Q$ is a finite dimensional
semi-simple \Q-algebra. For any $t\in \Z\langle S\rangle$ the
characteristic polynomial of $t$ on $H_i(X,r)\bigotimess L_r$  is independent
of $r$ and has rational coefficients by [loc.~cit.]. As $t$
stabilises a $\Z_{\ell}$-
(resp.~$\W(\k)$-) lattice in $H_i(X,r)\bigotimess L_r$, those coefficients are
$r$-integral for all $r$ and so integral. By the Cayley-Hamilton theorem the
image of $t$ in  $A\bigotimess \Q$ is integral over
\Z\ and so  $A$,
being equal to $\im:\Z\langle S\rangle\to A\bigotimes \Q$, is finitely
generated as it is contained in the different ideal of any order
containing it.
\end{lemma}

If still $L$ is algebraically closed of characteristic 0 and
$R$ and $T$ are two
sets, then we let $K(R,T,L)$ denote the Grothendieck group of the category of
finite dimensional $L$-representations of
$R \disjunion T$ such that every
element of $R$ commutes with every element of $T$. It is easy to see that every
simple object of this category is a tensor product of an irreducible
representation of $R$ and one of $T$ and that the two factors are
well-determined up to isomorphism. Hence $K(R,T,L)=K(R,L)\bigotimess K(T,L)$.
Furthermore, as $K(R,L)$ has a canonical base consisting of irreducible
representations we get a canonical pairing
$$
K(R,L)\bigotimess K(R,L)\to \Z
$$
where the the canonical base is orthonormal. Using this we get a mapping
$$
K(R,L)\bigotimess K(R,T,L)=K(R,L)\bigotimess K(R,L)\bigotimess K(T,L)\to K(T,L)
$$
and so for each $N\in K(R,T,L)$ a mapping
$$
N\cap  \co K(R,L)\to K(T,L).\tag1.4
$$
This mapping is compatible, in the obvious way, with the mappings obtained
from homomorphisms $L\to L'$ of algebraically closed fields. Hence we get
\begin{corollary}1.5.
Let $S'$ be a set of \k-correspondences of $X$ and suppose that
every element of $S$ commutes up to homological equivalence with $S'$. Then
\(1.4) gives a $Gal(\bQ/\Q)$-equivariant homomorphism
$$
[H^i(X)]\cap : K(S,\bQ)\to  K(S',\bQ)\tag(1.6)
$$
which equals, for each $r$, the restriction of
$[H^i(X,r)\bigotimess  L_r]\cap $ to K(S,$\bQ$).
\pro
\end{corollary}
\end{section}
\begin{section}Varieties of CM-type.

Let us fix an $n\in \N$ and assume, for simplicity, that
$$b_n(X)  =\sum _{i+j=n}\dim_kH^i(X,{\gOm}^j_{X/k}),$$
where $b_n(X) :=
\dim_{ L{_r}} H^n(X,r)$ for any $r$.
(This is of course always true
when $p=0$.) If $p>0$ this implies (cf.~\[Ek: IV, 1.2] or \[B-Og:\S8]) that
$H^n_{DR}(X/k)=H^n(X,p)/pH^n(X,p)$ and that if
$$
M^i:= \im(F^{-1}p^iH^n(X,p)\to  H^n(X,p)/pH^n(X,p))
$$
then
$$
M^i/M^{i+1}=H^{n-i}(X,{\gOm}^i_{X/k}).
$$
Hence, no matter the value of $p$, $H^n_{DR}(X/k)$
has a Hodge filtration with
the $H^{n-i}(X,{\gOm}^i_{X/k})$ as successive quotients.
\begin{definition}2.1.
$(X,S)$ is said to be of {\deffont separable CM-type in degree $n$} if the
$A_i$ of
\(1.3) has the property that $A\bigotimess \Z_{(p)}$ is a
separable (cf.~\[D-I:II,1])
$\Z_{(p)}$-algebra and for every $0\le i\ne j\le n$,
$H^{n-i}(X,{\gOm}^i_{X/k})$ and $H^{n-j}(X,{\gOm}^j_{X/k})$ are disjoint
$S$-modules (i.e.~they have no common composition factors).
\end{definition}
\begin{remark}
As the quotient of a separable algebra is separable it suffices to
verify that
$\Z_{(p)}\langle S\rangle$ factored by some known relations is separable.
\end{remark}
\begin{example}i) If $X$ is an abelian variety of CM-type in
the usual sense and
$End_k(X)$ is separable at $p$, which is always true if $p=0$, then
$(X,End_k(X))$ is of CM-type in degree 1 (and in fact in all other degrees).

ii) Kummer surfaces associated to abelian surfaces of CM-type are of
CM-type in degree 2. Hence, by \[S-I], K3-surfaces in characteristic 0 for
which the rank of the N{\'e}ron-Severi group is 20 are of CM-type in degree 2.

iii) (Fermat hyper-surfaces, diagonal automorphisms). This is well known
(cf.~e.g.~\[Ka:Sect.~6]).
\end{example}
\begin{lemma}2.2.
Suppose $(X,S)$ is of separable CM-type. If $p=0$ then $[H^n(X)]$ is
realisable over \Q\ and if $p>0$ then $A^i\bigotimess _\Z\Z_p$
is unramified (i.e.~a product of matrix algebras over unramified
extensions of $\Z_p$) and in particular $[H^n(X)]$
is realisable over $\Q_p$.
\pro The case $p=0$ follows by transcendental methods, in fact $[H_n(X)]$ is
realised by singular cohomology, and the $p>0$ is well-known (use
the fact that the Brauer group of a finite field is trivial and lift an
idempotent).
\end{lemma}
\begin{remark}
Is it possible to give an algebraic proof of the first part of the
lemma? The existence of \l-adic cohomology implies that it suffices to
prove realisability over \R.
\end{remark}
Suppose now that $(X,S)$ is of separable CM-type and let $M$ be an
irreducible component of $[H_n(X)]$. If $p=0$ there is then an irreducible
$S$-module $N$ such that $ M\bigotimess _\bQ\bk$ is a factor of
$M\bigotimess _\k \bk$ for an embedding of $\bQ$ in \bk\ and $N$ is a
sub-quotient of $H^n_{DR}(X/k)$.  (Note that the base extension of
$[H^n(X)]$ to \bk\ equals the extension of $[H^n_{DR}(X/k)]\in K(S,\k)$ to
\bk, which is seen by either using a constructibility argument to reduce to
\(1.3) or a transcendental argument.) By assumption there is a unique $i$,
$0\le i\le n$, such that $N$ occurs as a sub-quotient of \hodg n-i.i.  and
this $i$ depends only on $M$. If $p>0$ we get in the same way an
irreducible $K\langle S\rangle$-module $N$, $K:=\W(\k)\bigotimess \Q$, such
that $M\bigotimess _\bQ\bar K$ occurs in $N\bigotimess _{K}\bar K$ for an
algebraic closure $\bar K$ of $K$ and an embedding of \bQ\ in $\bar K$ and
$N$ occurs in $H^ n(X,p)\bigotimess _{\W(k)}K$. As $A\bigotimess \W(\k)$ is
separable, there is a unique, up to isomorphism, $A\bigotimess
\W(\k)$-lattice $N'$ with $N'\bigotimess _\W k=N$ (this follows from
\(2.2)) and $N'\bigotimess \k$ is an irreducible $A\bigotimess \k$-module
and we see that $N'\bigotimess \k$ is an irreducible sub-quotient of of
$H^n_{DR}(X/\k)$. By assumption there is then a unique $i$, $0\le i\le n$,
such that $N\bigotimess \k$ occurs in a \hodg n-i.i. and this $i$ depends
only on $M$. In both cases we put ${\gta}(M) := i$.

In conclusion we have obtained a mapping
$$
{\gta}: \Irr([H^n(X)])\to  n+1\;(:=\{0,1,\prickar ,n\}).
$$
Note that the action of $Gal$(\bk/\k) on $\ovl \Irr([H^n(X)]$) obtained through
the action of $Gal(\bQ/\Q)$ on it and the induced map $Gal(\bk/\k)\to
Gal(\bar\Q/\Q)$ (resp.~$Gal(\bk/k)\to Gal(\W(\bk)/\W(k))$)
preserves the fibers of {\gta} by construction.
\begin{definition}2.3.
Under the assumption of \(2.1) the type of $(X,S)$ is the
pair ($[H^n(X)]$,{\gta}).
\end{definition}

Finally, when $p>0$ the condition that $A$ is separable at $p$ implies that the
action of $Gal(\bar \Q/\Q)$ on $\Irr([H^n(X)])$ is unramified at $p$ so that we
may unambiguously speak about the action of the Frobenius morphism on
 $\Irr([H^n(X)])$ having once and for all chosen an embedding of  \bQ\
in $L_p$.
This permutation of $Irr([H^n(X)])$ we will denote {\gsi}.

We have now come to the main result of the present paper. Before we
formulate it we will need to introduce some constructions. Let $T$ be a set,
$M\in \ovl K(T,Q)$  and ${\gta}\co \Irr(M)\to  n+1$ a function.

Let us choose an embedding of \bQ\ into \C\ and let ${\gio}\in Gal(\bQ/\Q)$ be
the element corresponding through this embedding to complex
conjugation. Suppose
that for every ${\grh}\in \Irr(M)$, ${\gta}({\gio}({\grh}))=n-{\gta}({\grh})$
and also that $M$ is realisable over \Q\ by a module $V$.
For each simple factor
$A_r$ of $A_\Q$ we let $V_r$ be an irreducible $A_r$-module. We then put a
rational Hodge structure on $V_r$, of weight $n$, as follows: For ${\grh}\in
\Irr(M)\bigcap \Irr(V_r\bigotimess \C)$ we let $V_{r,{\grh}}$ be the
{\grh}-isotypical component of $V_r\bigotimess \C$ and then we put
$(V_r\bigotimess \C)^{i,n-i}:=\sum _{{\gta}({\grh})=i}V_{r,{\grh}}$.  We then
put a Hodge structure on $V$ by forcing $V=\bigoplus V_i\bigotimess
_{End(V_i)}Hom(V_i,V)$ to be an isomorphism of Hodge
structures. By construction
$T$ acts as morphisms of Hodge structures. However, the Hodge structure itself
depends only on the action of $Z(A_\Q(M))$ on $V$ so that an alternative method
of construction is to start with the set $\Irr(M)$ with its action of
$Gal(\bar\Q/\Q)$, let $K$ be the associated {\'e}tale
\Q-algebra, let $V$ be the
$K$-module of dimension specified by $M$ and then let $(V\bigotimess
\C)^{i,n-i}:=\sum _{{\gta}({\grh})=i}V_{\grh}$, where {\grh} runs over the
\Q-algebra homomorphisms $K\to \C$. In this way it is seen that $V$ as a
rational Hodge structure depends only on the
$Gal(\bQ/\Q)$-set $\Irr(M)$ and two
functions ${\gta}\co \Irr(M)\to n+1$ and $\dim\co \Irr(M)\to \N$, where $\dim$
is defined by $\dim(n)=dim_{Z(A_\C(n))}N$ for a representative of $N$ (\C\ can
of course be replaced by any algebraically closed field). If \k\ is a subfield
of \C\ such that the action of $Gal(\bk/\k)$ on $\Irr(M)$ preserves the fibers
of {\gta} then for a choice of descent of the Hodge filtration on each
$V_r\bigotimess \C$ to $V_r\bigotimess k$ for each we get a
descent of the Hodge
filtration on $V\bigotimess \C$ again by forcing the isomorphism above to
preserve the descent.

If $p$ is a prime such that $A(M)\bigotimess \Z_{(p)}$ is
finitely generated  as
$\Z_{(p)}$-module and separable
we associate, in a similar way, an $F$-crystal to $(M,{\gta})$:
Suppose that the action of $Gal(\bk/\k)$ on $\Irr(M)$ preserves
the fibers of {\gta}.  We know that $M$ is realisable
over $\Q_p$ and there is, up to isomorphism, a unique
$A(M)\bigotimess  \Z_p$-lattice $V$ such that $V\bigotimess \Q$ is such a
realisation. We also get analogous $V_r$. Further,
$V_r\bigotimess _{\Z_p}\W(k)$ is the sum of its isotypical components
$(V_r\bigotimess _{\Z_p}\W( k))_{\grh}$. The {\gsi}-linear
isomorphism $1\bigotimess {\gsi}$ takes
$(V_r\bigotimess _{\Z_p}\W(k))_{\grh}$ to $(V_r\bigotimess
_{\Z_p}\W( k))_{{\gsi}({\grh})}$
and we define the structure of an F-crystal on $V_r\bigotimess _{\Z_p}\W(k)$
by $F=p^{{\gta}(t)}(1\bigotimess {\gsi})\co (V_r\bigotimess
_{\Z_p}\W(k))_{\grh}\to (V_r\bigotimess _{\Z_p}\W(k))_{{\gsi}({\grh})}$.
The $F$-crystal structure on $V\bigotimess \W(k)$ is constructed as before.
Again $T$ acts by endomorphisms and there is an alternative way of
constructing the $F$-crystal if one is prepared to forget
the $T$-action. Indeed,
consider the set $\Irr(M)$ with its action of {\gsi} and the two functions
{\gta} and $\dim$.
We let $R$ be the set containing for each $n\in \Irr(M)\dim(n)$ copies of
$n$ with {\gsi} and {\gta} extended in the obvious way. We then consider
$\W(\bk)[R]$, the free $\W(\bk)$-module on $R$, and define the Frobenius map
by $F\vha r\hha=p^{{\gta}(r)}\vha {\gsi}(r)\hha$. Also if we
have for each $V_i$
an $End(V_i)$-representation {\grh} of $Gal(\bk/\k)$ on $V_i$, then we can
twist by this by letting $F$ act by $p^{{\gta}(t)}({\grh}\bigotimess {\gsi})$.

Finally, we would like to associate to $(M,{\gta})$ the \l-adic
analogue of this, that is a Hecke character.
We will be able to associate to our data the algebraic part of a Hecke
character but a problem arises as there is no canonical choice for a Hecke
character with a given algebraic part, indeed such a character may exist
only after an extension of the coefficient field. This will have as a
consequence that our description of the \l-adic cohomology will not be as
satisfactory as the description of the Hodge structure or $F$-crystal of a
variety of CM-type. In case $[H_n(X)]$ is multiplicity free we will be able
to do better however.
In any case the algebraic part (cf.~\[De:5.3]) can be associated
to our data as follows. Assume that for any embedding of \bQ\ in
\C\ we have ${\gta}({\gio}({\grh}))=n-{\gta}({\grh})$
as above for the corresponding {\gio}. Assume also that \k\ is a number field
for which the action of $Gal(\bQ/K)$ on $\Irr(M)$ stabilises
the fibers of {\gta}.
 If again  $K$ is the \Q-algebra
corresponding to  $\Irr(M)$, then we can define a multiplicative map
$$
\eqalign{
k^\times &\;\longrightarrow\;\; K^\times  \cr
{\gla}&\mapsto \prod _{{\grh}\in \Irr(M)}{\grh}(N_{k/\Q}({\gla}))
^{{\gta}({\grh})}.\cr}
$$
By the assumption on \k\ this is well-defined and by the assumption on
{\gta} the projection onto each simple factor of $K$
fulfills the conditions for being the algebraic part
of a Hecke character of weight $n$ so we obtain in this way a set of algebraic
parts of Hecke characters.
\begin{theorem}2.4. Let $(X,S)$ be of separable CM-type in degree $n$.

i)  If $k\subseteq \C$
then $H^n_{sing}(X(k),\Q)$ is isomorphic as a Hodge structure
with $S$-action to the one associated to the type of $(X,S)$ with a descent of
the Hodge filtration to \k\ of the sort described.

ii)   If $p>0$ then $H^n(X,p)$ is isomorphic as $F$-crystal to the $F$-crystal
associated to the type of $(X,S)$ and a representation of $Gal(\bk/\k)$.

iii)  After a finite extension of \k\ the
$Gal(\bk/\k)$-representation on $H_n(X,r)$, $(r\ne p)$
factors through the Galois group of the algebraic closure
$K$ of the prime field
in (the finite extension of) \k. If $p=0$ this representation is given, on
the Galois group of a finite extension of $K$, by a direct sum of
algebraic Hecke characters whose algebraic
parts are the ones associated to the type of $(X,S)$ and with multiplicities
given by $\dim$.

iv) If $\Irr(M)$ is multiplicity free (i.e.~every
irreducible representation of $S$
occurs at most once in $\Irr(M)$) and \k\ is a number field then the
$Gal$(\bk/\k)-representation on $H_n(X,r)$ is given by a direct sum of
algebraic Hecke characters with values in the simple components of
$Z(A_\Q(\Irr(M)))$.

\pro To begin with let $R$ be \Q, $\Z_p$ resp.~$\Q_r$. Then $S$ generates an
$R$-subalgebra $B$ of the algebra of endomorphisms of $H^n_{sing}(X(\k),\Q)$,
$H^n(X,p)$ resp.~$H^n(X,r)\bigotimess _{\Z_r}\Q_r$ such that
$$
B/(\hbox{\rm maximal
nilpotent ideal})=A^n\bigotimess _\Z R
$$
(the $A^n$ being that of \(2.1)). By \[C-R:Thm. 72.19]
$B\to A^n\bigotimess  R$ splits as an algebra map and so we can assume that
$A\bigotimess  R$ acts on $H^n_{sing}(X(\k),\Q)$, $H^n(X,p)$ resp.
$H^n(X,r)\bigotimess _{\Z_r}\Q_r$. Using \(1.1) we reduce to the case when
$A\bigotimess R$ is a division algebra or, in the case of ii), isomorphic to
$W(\F)$ for a finite field \F. The proof of i) is then easy: By the comparison
theorem $H^n_{sing}(X(\k),\Q)$ is a representative of $[H^n(X)]$ and as
the Hodge decomposition on $H^n_{sing}(X(\k),\Q)\bigotimess _\Q\C$ is
stable under $A\bigotimess  \C$, the assumption of CM-type forces the Hodge
decomposition to be obtained as the lumping together of isotypical
components.

As for ii), the fact that $A\bigotimess \W(k)$ is separable
implies that we have a
unique isotypical decomposition $H^n(X,p)=\bigoplus  M_{\grh}$, where {\grh}\
runs over the irreducible $A\bigotimess  K$-modules
($K:=\W(\k)\bigotimess _\Z\Q$)
occurring in $H^n(X,p)\bigotimess \Q$. As $F$ is \gSi-linear
and commutes with $A$ it
maps $M_{\grh}$
to $M_{\gSi({\grh})}$. If $(p^{n_1},p^{n_2},\ldots,p^{n_k})$
are the elementary divisors of the linear mapping $F\co
M_{\grh}\to\gSi_* M_{\gSi({\grh}))}$, the characterisation of the Hodge
filtration recalled at the beginning of this section shows that
$M_{\grh}/pM_{\grh}$ is non-disjoint from \hod
X{n-i}ik exactly when $i$ equals
some $n_j$. Hence, by assumption, all the $n_j$ equal ${\gta}({\gLa})$ where
{\gLa} is any component of $[H_n(X)]$ which occurs in
${\grh}\bigotimess _K\bar K$.
We will denote, by abuse, this common value ${\gta}({\grh})$. Thus $F\co
M_{\grh}\to\gSi_* M_{\gSi({\grh})}$ is p$^{{\gta}({\gLa})}$ times an
isomorphism. Hence if we define $F'\co H^n(X,p)\to H^n(X,p)$ as being
$p^{-{\gta}({\grh})}F$ on $M_{\grh}$, $H_n(X,p)$ becomes a unit root crystal
and is hence described by the $Gal$(\bk/\k)-representation on the fixed
points of $F$ (over \bk). This action commutes with
$A\bigotimess _\Z{\Z_p}$ and
so gives the desired description.

Let us now turn to iii) and let us begin with the case $p=0$. We may assume
that \k\ is a finitely generated field. We have a representation ${\gph}\co
Gal(\bar k/k)\to Aut_{Z\bigotimess \Q_r}(H^n(X,r)\bigotimess
\Q)$ and after possibly
enlarging $Z$ ($:=Z(\Irr(M))$) and \k\ we may assume that there exists an
algebraic Hecke character $I_m(k')\to Z^\times $, where $k'$ is the
algebraic closure of \Q\ in \k, whose algebraic part is the one coming from
the type of $(X,S)$ (using that the condition on {\gta} is fulfilled by the
transcendental theory). Twist {\gph} by this character, considered as a
character of $Gal(\bk/\k)$ through the
morphism $Gal(\bk/\k)\to Gal(\bk'/\k')$ and (\[Se:II,2.7]). What we now need
to prove is that this twist ${\gph}'$ has finite image (cf.~\[De:Thm.~5.10]).

Let us first show that if $\gSi\in Gal(\bar k/k)$ then ${\gph}'(\gSi)$ is
quasi-unipotent. The possible orders for the eigenvalues of a
quasi-unipotent matrix over $\Q_r$ of given order is bounded as the
degrees of the extensions of $\Q_r$ obtained by adjoining an $m$th root of
unity goes to infinity with $m$. It is
therefore enough to verify the quasi-unipotence on
a dense set of Frobenius elements. By the \v Ceboratev density theorem it
suffices to check quasi-unipotence for the Frobenius elements corresponding
to maximal ideals for some thickening of \Sp\k\ over which $A_n$ is
separable, $X$ is smooth and
$$
b_n(X)=\sum _{i+j=n}dim_h\hod Xijh,
$$
where $h$ is the residue field. If $F_m$ is the Frobenius element of
$Gal(\bk/\k)$ we then want to show that all the eigenvalues of
${\gph}'(F_m)$ are roots of unity or, as they are all algebraic numbers, that
all their absolute values are equal to 1. For the infinite primes we use the
Riemann hypothesis for $X$. At finite places away from $q:=\symb{char} h$
there is no
problem. Let us therefore consider the places over $q$. Pick
a place $v$ of $Z$ lying
over $q$ normalised so that $v(|h|)=1$. By definition $v({\gph}(F_ m))$ equals
the average of {\gta} over the orbit of $v$ of the action of \gSi\ on
$\Irr([H^n(X)])$, where $v$ is seen as a homomorphism $Z\to L_q$ and
thus giving an element of $\Irr([H^n(X)])$ (recall that $L_q$ is an algebraic
closure of $\Q_q$). Hence we want to show that any eigenvalue of the
action of $F_m$ on the $v$-isotypical part of $H^n(X,q)$ has the same
valuation. By construction $(X_h,S)$ is of separable CM-type in degree $n$ and
the eigenvalues of $F_m$ are of course the same for $X_k$ and
$X_h$ so we may replace $k$ by $h$. Applying \(1.5) to $S$
and $\{F_h\}$ we see, as
$F_m=F^*_h$ on $H^n(X,r)$, that the eigenvalues of $F_m$ on the
$v$-isotypical part of of $H^n(X,r)$ are the same as the the eigenvalues of
$F^*_h$ on $H^n(X,q)$. By ii), if $u$ is the length of the \gSi-orbit of $v$
then $F^u$ is divisible exactly by $p^t$ on the $v$-isotypical component of
$H^n(X,q)$, where $t$ is the sum of the values of {\gta} over the \gSi-orbit
of $v$. As $F^*_h=F^r$, where $|h|=p^r$, we immediately get what
we want.

Now again as the orders of the eigenvalues of the elements of
$Gal$(\bk/\k) are bounded after replacing \k\ by a finite extension we
may assume that the image of ${\gph}'$ consists entirely of unipotent
matrices and so by Engel's theorem ${\gph}'$ is a unipotent representation.
We aim to show that it is in fact trivial. As
$G:={\gph}'(Gal(\bar k/k))$ is a compact $r$-adic Lie group the closed subgroup
of $G$ generated by $r$th powers is of finite index in $G$ and by the Frattini
lemma any closed subgroup mapping surjectively onto the quotient of $G$ by
this subgroup equals all of G. We may then apply the Hilbert irreducibility
theorem  to get a number field specialisation
$\k''$ of \k\ such that $X$ has good reduction at $\k''$ and that the composed
map $Gal(\bk''/\k'')\to Gal(\bk/\k)\to G$ is surjective. Hence we may
assume that \k\ is a number field. I claim that for each prime of \k\ over
$r$, the inertia group of that prime has finite image in
$G$. Indeed, ${\gph}'$ is
Hodge-Tate as an algebraic Hecke character is always Hodge-Tate
and by \[Fa]. The finiteness then follows from
(\[Se1:1.4,Cor.~3]) as the unipotence
implies that the Hodge-Tate weight is zero. As ${\gph}'$ is unipotent this
implies that ${\gph}'$ is unramified over $r$. The other monodromy groups
automatically have finite, and therefore trivial, images. The triviality of $G$
then follows from the finiteness of the Hilbert class field of \k.

We have therefore proved iii) when $p=0$. The case $p>0$ is similar up to
the point where we have arrived at a unipotent representation.  Any
homomorphism $Gal(S,\bar s)\to \Z_r$ for a finitely generated
$\F_p$-scheme $S$ is geometrically trivial, by
\[K-L:Thm~1], so by thickening \k\ we finish.

As for iv) we start as above so that we have an action of
$A:=A^n\bigotimess  \Q$ on $H^n(X)$ by correspondences. Note
that the assumption of multiplicity freeness implies
that the commutant of $A$
in $End(H^n(X,r))$ equals $Z\bigotimess  \Q_r$. Let $v$ be a
place of \k\ at which
$X$ has good reduction with fiber $X_v$ over the residue field $\F_v$. Apply
the construction of a semi-simple algebra of correspondences to all
correspondences so as to get $B$. Then $B$ contains the Frobenius
correspondence in its center as well as the subalgebra  $A$. Let $C$ be the
commutant of $A$ in $B$ so that $B=AC$. By the observation
just made $C\bigotimess
\Q_r=Z\bigotimess  \Q_r\subseteq  B\bigotimess \Q_r$ and so $C\subseteq B$ and
therefore $A=B$. Hence the Frobenius correspondence $F_v$ lies in $Z$.
Therefore we have associated to every place $v$ of \k\ outside a finite set an
element $F_v$ of $Z$. Extending by multiplicativity we get a homomorphism
$I_m(k)\to Z^\times$ for a suitable $m$. As in the proof of iii) we show that
the projections onto the simple factors of $Z$
are algebraic Hecke characters with algebraic parts given by the type of
$(X,S)$.
\end{theorem}
\begin{remark}
i) It is probably true that in ii) we also get the conclusion
that $[H^n(X,p)]$ is geometrically constant. This would
follow from a good theory of
over-convergent $F$-crystals.

ii) Can one find a good extension of iii) that would contain iv) as a special
case?

iii) An example showing that there are problems in the \l-adic case is
obtained as follows. Pick an imaginary quadratic field $K$ with class number
greater than one. There is an elliptic curve $E$ with complex multiplication by
the ring of integers $R$ of $K$ defined over the Hilbert
class field $H$ of K. The pair
($R_{H/K}E$,R) is of CM-type in degree 1 over $K$ yet there is no algebraic
Hecke character whose algebraic part is that obtained from the type of
($R_{H/K}E$,R).
\end{remark}
As will come as no particular surprise, for abelian varieties our notion
coincides with the traditional one.
\begin{proposition}2.4.
Suppose $(X,S)$ is of separable CM-type in degree 1. Then its
Albanese variety is of CM-type in the usual sense possibly after a finite
extension of \k.

\pro This follows from \(2.4) and the Tate conjecture for homomorphisms
between abelian varieties. Another proof is for $p=0$ to
note that \(2.3:i) says
that the degree 1 Hodge structure of $X$ is visibly of CM-type and for $p>0$
that $(Alb X,S)$ is rigid as by definition $Hom_S(\hod
X01k,\hod X10k)$ is equal
to 0 and so after a finite extension of \k, $Alb\, X$ can be defined over a
finite field and is hence of CM-type.
\end{proposition}
\end{section}
\begin{section}Hereditary CM-type

\mdef\End{\symb{End}\,}
Theorem 2.3 suffers somewhat on the $p$-adic side as the very natural
example of $(E,\End E)$ where $E$ is a supersingular elliptic curve is not
of separable CM-type; $\End E$ is not separable at $p$. It is possible to
give a result which in that case specialises to a satisfactory answer. In
this section we will give a generalisation of the previous results that
will cover this case. The maximal possible generality would seem to be to
assume that $A^n\bigotimess \Z_{(p)}$ should be a {\deffont hereditary}
order which means that any $A^n$-splitting of crystalline cohomology
tensored with \Q\ comes from an $A^n$-splitting of crystalline cohomology
itself. Let us recall that an order is hereditary if each lattice over it
is projective.
\begin{remark}
The meaning of the term differs somewhat in various areas of the literature
as hereditary sometimes means just that a submodule of a projective module
is projective. The definition used here means that the base extension of
the order to the fraction field of its base ring is semi-simple together
with the fact that every sub-module of a projective module is
projective. We will want this extra condition and hence adopt the current
definition (which is to be found for instance in \[Re]).
\end{remark}
On the other hand, the example of an automorphism of order $p$ acting
(non-trivially) on a curve of genus $(p-1)/2$ shows that the condition that
different Hodge pieces be disjoint is not reasonable as the cyclic group of
order $p$ has only one irreducible representation mod $p$.  The situation will
no longer be as simple as in the separable case. It is still true that one to
any irreducible $A^n\bigotimess \k$-module can associate an irreducible
$A^n\bigotimess K$-module but this map is no longer injective (though
surjective). We will use \[Re:Ch.~9] as a general reference to the theory of
hereditary orders. For the reader's convenience we repeat the salient facts in
the following proposition as well as adding a result -- a weak version of the
elementary divisor theorem -- which is not to be found in
\[loc.~cit.] (but no doubt is
not new).
\begin{proposition}eldiv
Let $A$ be a hereditary order over a henselian discrete valuation ring $R$
with fraction field $K$.

i) Any submodule of an $A$-lattice is a submodule of finite colength of a
direct factor of the lattice.

ii) In every indecomposable $A$-lattice there is exactly one submodule of a
given colength.

iii) If $M$ is an indecomposable $A$-lattice and $M\hookrightarrow N_i$ two
inclusions of finite colength. Then one of these inclusions is contained in
the other.

iv) Every indecomposable finitely generated torsion $A$-module is a quotient
of an indecomposable $A$-lattice.

v) Let $M$ be an $A$-lattice and $N$ a sub-lattice of it. Then there is a
decomposition of $M$ as a direct sum of indecomposable submodules whose
intersection with $N$ also gives a decomposition of $N$ into a direct sum of
indecomposable submodules.
\pro For i) we take the saturation of the submodule. The quotient of the
lattice by that saturation is torsion-free and hence projective and the
saturation is therefore a direct factor. For ii) we notice that by i) any
submodule of the lattice is also indecomposable so we may assume by
induction that the given colength is 1. However, the lattice being
projective is the projective hull of its co-socle (the maximal semi-simple
quotient) and so being indecomposable the co-socle is simple which means
that there is a unique submodule of colength 1; the radical. For iii) we
note that $M\hookrightarrow N_i$ are included in a common inclusion of
finite colength (being of finite colength). We then apply ii).

As for iv) we use induction on the length of the module $M$. We therefore
find a simple quotient $S$ of $M$ and apply the induction hypothesis to the
kernel $M'$ of this map. We will temporarily (and improperly) call a
torsion quotient of an indecomposable lattice a {\it cyclic} module. Thus
we may assume that $M$ is an extension of a sum of cyclic modules by the
simple module $S$. This extension is the sum, as extension, of the
extension of the cyclic summands by $S$. Let us first study the latter
extensions and let us denote by $P$ the projective hull of $S$, by $Q$ its
radical, by $V$ the cyclic summand and by $R$ its projective hull. Then
every extension of $V$ by $S$ comes from pushout by a map from $Q$ to $V$,
the same extensions being obtained if the difference of two morphism
extends to a map from $P$ to $S$. Now, I claim that all non-surjective maps
$Q\to S$ so extend. In fact the map lifts to a map $Q\to R$ which
necessarily is injective as $Q$ is indecomposable. If the map $Q\to S$ is
not surjective then the map $Q\to R$ is neither. By applying iii) we see that
$P$ must be isomorphic to the unique submodule of $R$ containing $Q$ as a
submodule of colength 1 and thus the original map lifts to $P$. This result
shows that if $V$ is not a quotient of $Q$ then any extension of $V$ by $S$
is trivial and if it is, then the group of extensions can be identified
with maps from the co-socle of $Q$ to the co-socle of $V$ which are
isomorphic simple modules. We will now show that, after possibly changing
the direct sum decomposition of $M'$ we may assume that all but one of the
extension classes of direct summands by $S$ are trivial. This will clearly
show iv). For this we may immediately discard summands of $M'$ which are
not quotients of $Q$ as their extension classes have just been shown to be
trivial. Furthermore, we may use induction on the number of non-trivial
extension classes. Note now that if $V_1$ and $V_2$ are summands then for
any map \pil\phi{V_1}{V_2} we may consider the automorphism of $V$ which
maps $v\in V_1$ to $v+\phi v$ and acts as the identity on all other
factors. If $e_i$ are the extension classes then all of them but $e_2$ are
unchanged and $e_2$ is changed into $e_2+\phi_*e_1$. If we identify
extension classes of $V_i$ with homomorphisms from the co-socle of $Q$ to
that of $V_i$ then $\phi_*$ is just composition by the map on co-socles
induced by $\phi$.  As $V_1$ and $V_2$ are both both quotients of $Q$, by
ii) on is a quotient by the other and we may assume that $V_2$ is a
quotient of $V_1$. In that case, any endomorphism of $Q$ induces a morphism
$V_1\to V_2$ and by the projetivity of $Q$, any map from the co-socle of
$V_1$ to that of $V_2$ is induced by an endomorphism of $Q$. Putting this
together we see that any extension class is of the form $\phi_*e_1$ so that
the we may choose $\phi$ so that $\phi_*e_1=-e_2$ which allows us to
decrease the number of non-zero extension classes.

To finally prove v) we consider the module $M/N$. This is a direct sum of
a lattice and a torsion module and the lattice may be split off from $M$
without changing $N$. Thus we may assume that $M/N$ is torsion. We then use
iv) to write that quotient as a direct sum of quotients of indecomposable
projective modules. The sum of the projective hulls of each summand is a
projective hull of the sum. That projective hull is a direct summand of the
map $M\to M/N$. This immediately gives the pair $(M,N)$ as a direct sum of
of pairs $(P_i,P'_i)$, where $P_i$ is indecomposable and a factor
$(M',M')$. As $M'$ is a sum of indecomposables, v) follows.
\end{proposition}
We will need a definition which is very special to the situation at hand.
\begin{definition}
Let $A$ be a hereditary order over a henselian mixed characteristic discrete
valuation ring $R$ with positive residue field
characteristic $p$ and let $M$ be
a finitely generated torsion module killed by $p$. By the
\definition{complementary module} to $M$ we mean the torsion module (defined up
to isomorphism) obtained as $P/pP'$, where $P$ is a projective hull of $M$ and
$P'$ is the kernel of the natural map $P\to M$.
\end{definition}
As there can be, as opposed to the separable case, non-trivial extensions of
modules we also will need to recall the definition of block, well-known in the
theory of general orders,
\begin{definition}
Let $A$ be a hereditary order over a henselian discrete valuation ring $R$
with fraction field $K$. Two indecomposable (finitely generated) $A$-modules
belong to the same block if there is a non-zero morphism from the projective
hull of one to the other. This is equivalent to the two hulls tensored with $K$
being isomorphic. If $M$ is a finitely generated $A$-module then the
\definition{$B$-component}of $M$ is the sum of all indecomposable
factors belonging to the block $B$. (It is clear that any finitely generated
$A$-module is the direct sum of its components associated to different blocks.)
\end{definition}
What is different with the hereditary case as opposed to the separable case is
that we may have non-semisimple (f.g.) modules killed by $p$. This will imply
that to define CM-type it is not enough to look at what simple modules occur in
which Hodge piece; the more precise module structure needs to be taken into
account. As we will see this forces certain relations between Hodge pieces.
Our results will be purely algebraic so we will, rather than sticking
to the notation of this article as a whole, use the
following notation: $A$ will be a
hereditary $\Z_p$-order and $M$ will be an $F$-crystal with
an action of $A$. We
define the Hodge filtration on $M/pM$ by $M^i:=F^{-1}p^iM/pM$ and the Hodge
modules $H^i:=M^i/M^{i+1}$ (which may be considered as
$A\bigotimes \W$-modules).
\begin{definition-lemma}
We define the \definition{$A$-primitive part} of $H^i$ as the direct factor
(defined up to isomorphism only) by induction on $i$. For $i=0$ we let the
primitive part be all of $H^i$. For $i>0$ the complement of the primitive part
of $H^{i-1}$ is a direct factor of $H^i$ and we let the primitive part be a
complementary factor of it.
\pro What is to be proven is the statement about the
complement of the primitive
part being a direct summand. We will give another description of the primitive
part which will make this obvious. Consider therefore the Frobenius map as a
$\W$-linear map $\sigma^*M\to M$, which then also is a
$A\bigotimess \W$-linear. This is a hereditary order so we may by \(eldiv:v)
split this map up in indecomposable factors. Using Mazur-Ogus' characterisation
(\[B-Og]) of the Hodge filtration and the fact that submodules are linearly
ordered we immediately see that each indecomposable factor will contribute a
cyclic module to one Hodge piece and its complement to the next.
\end{definition-lemma}
 We are now ready to define what we
mean by CM-type in the context of actions of hereditary orders.
\begin{definition}
The pair $(M,A)$ is of \definition{hereditary CM-type} if $A\bigotimess
\Z_p$ is a hereditary order and for each block $B$, the
$B$-component of the primitive part of $H^i$ is non-zero for at
most one $i$ and all indecomposable factors of that $B$-component have the same
length.
\end{definition}
We have now set up our definitions so that we may carry through the same
analysis as in the separable case (it should be noted that in the case that
$A^n\bigotimess \Z_{(p)}$ is actually separable this definition coincides with
the previous one).
\begin{lemma}galois
Let $B$ be a $\W$-algebra, finitely generated and free as a
$\W$-module. Suppose
$n$ is a positive integer and $T$ a $\sigma^n$-linear
automorphism of $B$. Using
$T$ to get a \Z-action on $B$ we have that $H^1(\Z,B^\times)=*$.
\pro If $B'$ is the $\Z_p$-algebra of $T$-fixpoints then we have that
$B=B'\bigotimes \W$ and an element of  $H^1(\Z,B^*)$ is
given by an automorphism
class of a finitely generated right $B'$-module whose extension of scalars to
$\W$ is isomorphic as $B$-module to $B$ itself. As the
extension $\Z_p\to \W$ is
faithfully flat this means that such a $B'$-module is projective. Hence it is
determined up to isomorphism by its co-socle and to prove
the lemma it is enough
to show that if we have two semi-simple $B'$-modules which become isomorphic
under extension of scalars to $\W$ are isomorphic. This
however is obvious (using
for instance the independence of central characters of a semi-simple algebra).
\end{lemma}
\begin{theorem}
Suppose $(M,A)$ is of hereditary CM-type and that {\k} is algebraically
closed. Then it is determined up to isomorphism by which blocks appear in the
primitive part of which Hodge-modules and the common length of indecomposable
factors of each such block.
\pro Note first that we can make {\gsi} act on the blocks of
$A\bigotimes \W$ by
the condition that $N$ belongs to the block $B$ iff $\sigma^*N$ belongs to
$\sigma^*B$. If we now split up $M$ in blocks, $M=\bigoplus_B M_B$, then it
is clear that $F$, considered as a map $\sigma^*M\to M$ is a sum of maps
$\sigma^*M_B\to M_{\sigma^*B}$. Now, for an indecomposable $B\bigotimes
W$-lattice  $N$ the length of $N/pN$ only depends on which block $N$ belongs
to. Indeed, any two indecomposable $B\bigotimes W$-lattices in the same block
are contained in each other with quotient of finite length, the kernel and
cokernel of the map induced by reduction modulo $p$ then has the same length.
We now consider the component of $F$, $\sigma^*M_B\to M_{\sigma^*B}$, as a
$W$-linear map and split it up into indecomposable pieces according to lemma
\ref{eldiv}. Looking at each indecomposable piece we see that if $B$ appears in
the $i$'th Hodge piece of $M$ then $F$ maps $\sigma^*M_B$ into
$p^iM_{\sigma^*B}$, the image contains $p^{i+1}M_{\sigma^*B}$ and the length of
each indecomposable factor of $p^iM_B/M'$, $M'$ being the image, has the same
length (as each such length added to the common length of
the indecomposables of
the primitive $B$-part of the Hodge piece adds up to the common length of an
indecomposable lattice in $B$ modulo $p$). This means that any indecomposable
factor of $M_{\sigma^*B}/M'$ has the same length, which is the same as saying
that $M'=\symb{rad}^mM_{\sigma^*B}$ for a suitable $m$,
where $\symb{rad}(-)$ is
the radical functor. We may thus use $F$ to identify $\sigma^*M_B$ with
$\symb{rad}M_{\sigma^*B}$, where $m$ is determined by $i$ and the common length
of indecomposables of the primitive part belonging to the block $B$. If $n$ is
the smallest positive integer for which $\sigma^{n*}B=B$, then $F^n$ maps $M_B$
onto $\symb{rad}^kM_B$ for a suitable $k$. It is then enough to show that
all such maps are conjugate under automorphisms of $M_B$. Fix one such map
{\gph}. Now, I claim that the relation $\phi\circ f^\sigma=g\circ\phi$ defines
an automorphism $f\mapsto g$ of $\End(M_{\sigma^*B})$.
Indeed, for any $g$ there
is an $f$ fulfilling that relation as the image of {\gph} is equal to
$\symb{rad}^mM_{\sigma^*B}$. Conversely, the inverse image of $M_{\sigma^*B}$
in $\sigma^*M_B\bigotimes K$ under {\gph} is equal to its
largest sub-lattice for
which the quotient by $\sigma^*M_B$ has all its indecomposable components of
length less than or equal to $m$ which shows that to any $f$ there is a $g$. As
any map fulfilling the conditions imposed on {\gph} differs from it by an
automorphism of $\sigma^*M_B$ we can apply lemma \ref{galois} to include that
there is, up to isomorphism, only one $F$.
\end{theorem}
I would also like to record that, just as in the separable case, multiplicity
freeness implies CM-type.
\begin{proposition}
Suppose that $(M,A)$ is multiplicity free in the sense that an irreducible
$A\bigotimes K$-module appears at most once in $M\bigotimes K$. Then $(M,A)$ is
of hereditary CM-type.
\pro The condition implies that for any block, the component of $M$ in that
block is indecomposable. That immediately implies that a given block appears in
the primitive part of just a single Hodge piece and that part is indecomposable
so the condition on length is fulfilled.
\end{proposition}
We finish this section with some examples.
\begin{example}
i) Consider a supersingular elliptic curve $E$ and its ring
of endomorphisms $A$,
which is a maximal order in a division ring and hence hereditary. The
action of $A$ on the first crystalline cohomology group is
multiplicity free and
hence of hereditary CM-type. More precisely, $H^1(E,p)$ is an indecomposable
$A\bigotimes W$-lattice and $\sigma^*H^1(E,p)$ is the other indecomposable
$A\bigotimes W$-lattice -- both are in the same block. Hence, $H^0(E,\Omega^1)$
is one irreducible $A\bigotimes W$-module and $H^1(E,\Cal O_E)$ the other. The
image of $\sigma^*H^1(E,p)$ under $F$ is the maximal proper submodule of
$H^1(E,p)$.

ii) Let $C$ be the projective, smooth completion of the curve $y^p-y=x^2$,
$p\ne2$, and consider the action of $\Z/p$ given by $y\mapsto y+\alpha$. The
action of the group algebra of $\Z/p$ on $H^1(C,p)$ factors
through the quotient
$A$ that is the ring of $p$'th roots of unity. $H^1(C,p)$ is then a free
$A\bigotimes \W$-module of rank 1 and hence $(C,A)$ is of hereditary
CM-type. This time the situation is simpler as the ring is commutative and to
prove the classification theorem we could simply have divided $F$ by
$(\zeta-1)^{(p-1)/2}$ to obtain a unit root crystal.
\end{example}
\end{section}
\begin{section}Applications to the N{\'e}ron-Severi group.

In this section we will suppose that $X$ is a surface and that $(X,S)$ is of
separable type in degree 2. Then $S$ acts on the N{\'e}ron-Severi group $NS$ of
$X_{\bar k}$. If $p=0$, \(2.4) and the Lefschetz theorem on $(1,1)$-classes
show that $[NS\bigotimess \bar\Q]\in K(S,\bQ)$ equals the sum of all
irreducible {\grh}\ in $[H^2(X)]$ such that ${\gta}({\gsi}({\grh}))=1$ for all
${\gsi}\in Gal(\bar\Q/\Q)$ and the Tate conjectures implies this in all
characteristics. However, \(2.3) can be used to obtain further
information on $NS$. To illustrate this let us suppose that $p>0$ and that
$\rk NS=b_2(X)$. By possibly extending \k\ we may assume that $NS$ is
defined over \k. As $c_1\co NS\bigotimess \Z_\ell\to H^2(X_{\bk},\Z_\ell)$
($\ell\ne p$) has torsion free cokernel (cf.~\[Gro:8.7]) and
the two modules have
the same rank, $c_1$ is an isomorphism. By Poincar{\'e} duality the
intersection pairing is perfect at \l. By \[Ill:II,5.8.5,5.20] the image of
$c_1\co NS\bigotimess  \W(k)\to H^2(X,p)$ is the largest
sub-$F$-crystal in which $F$ is
divisible by $p$ and, again by Poincar{\'e} duality, if ${\gsi}_0$ is the
$\W(\k)$-length of the cokernel, then $p^{2{\gsi}_0}$ is the exact power of $p$
dividing $disc(NS)$. Hence  by the Hodge index theorem
$disc(NS)=(-1)^{b_2-1}p^{2{\gsi}_0}$.

As the whole $H^2(X,p)$ is determined by the type of $(X,S)$, ${\gsi}_0$ is as
well and we will now see how this can be done explicitly. By \(2.3)
$M:=H^2(X,p)=\bigoplus _{{\grh}\in \Irr([H^2(X)])}M_{\grh}$
and $F\co M_{\grh}\to
M_{{\gsi}({\grh})}$ is $p^{{\gta}({\grh})}$ times an isomorphism.
Let $N\subseteq M$ be the maximal sub-$F$-crystal on which $F$
is divisible by $p$. Consider
 $T:=\Irr([H^2(X)])$ with the functions {\gta} and $\dim$ and the action of
{\gsi}. We shall now describe an algorithm for computing
${\gsi}_0$. To do this we
start by by considering $M$ with $F':=p^{-1}F$ as a virtual $F$-crystal
i.e.~$p^{-1}F$ takes $M$
into $M\bigotimess \Q$ rather than into $M$ itself. Now $N$ can then be
characterised as the maximal sub-$F$-virtual-crystal which is actually a
crystal. As it is unique it is a sub-representation and so it is the direct
sum of the $N_{\grh}$. We will now concentrate on one specific {\gsi}-orbit
on $T$ and assume that $M$ is in fact the $F$-crystal associated to
it. Pick one {\grh} in this orbit. As $M_{\grh}$ is of rank 1 $N_{\grh}$ is
equal to $p_nM_{\grh}$ for some $n$. All powers of $F'$ must take
$N_{\grh}$ to $M$ which means that $n+\sum_{j=0}^k({\gta}({\gsi}^j{\grh})-1)$
is greater than or equal to 0 for all $k$. Hence if we put $n$ equal to
$-\min_k\sum _{j=0}^k({\gta}({\gsi}^j{\grh})-1)$ and define $N'$ as
$\bigoplus p^{n_{{\grh}'}}M$, where $m_{{\grh}'}:=n+\sum
_{j=0}^k({\gta}({\gsi}^j{\grh})-1)$ with ${\grh}'={\gsi}^k{\grh}$ we have a
sub-$F$-virtual-crystal of $M$ which clearly is an actual $F$-crystal (here we
use the fact that the sum of ${\gta}-1$ over the orbit is 0). We have
also seen that $N_{\grh}\subseteq N'_{\grh}$ and as $N$ is the maximal
sub-$F$-crystal we have equality.  Finally, again using that the sum of
${\gta}-1$ over the orbit is 0 it is immediately realised that $N'$ is
independent of the choice {\grh} and so has to be equal to $N$. In
particular we see that the contribution of this orbit to ${\gsi}_0$ equals the
multiplicity of the orbit times the sum of the $m_{{\grh}'}$.

\begin{example}
We consider one orbit for {\gsi} and describe such an orbit by
$({\gta}(t),{\gta}({\gsi}(t)),\ldots,{\gta}({\gsi}^{h-1}(t)))$,
where $h$ is the
length of the orbit. We also assume that the starting point {\grh} is the first
element of this list.

\noindent i)\ (0,2) gives partial sums $(-1,0)$ and so $n=1$
and the list of the
$m_{{\grh}'}$ is $(0,1)$ and finally the contribution to ${\gsi}_0$ is 1.

\noindent ii)\  $(0,2,1)$ gives partial sums $(-1,0,0)$, $m$s $(0,1,1)$ and
${\gsi}_0=2$.

\noindent iii)\ $(2,1,1,1,1,0)$ gives partial sums $(1,1,1,1,1,0)$, $m$s
$(1,1,1,1,1,0)$ and so ${\gsi}_0=5$.
\end{example}
    As a geometric example let us first consider the Fermat surface $X_m=
\{X^m_0+X^m_1+X^m_2+X^m_3=0\}$ and the group of diagonal automorphisms
$G_m=\mu^4_m/(scalars)$.
The irreducible representations of this group are the
elements of the dual group
$$
\check G_m:=\{(b_0,b_1,b_2,b_3)\in(\Z/m)^4:\sum _{i=0}^3b_i=0\}
$$
and it is well known (cf.~\[Ka:Sect.~6]) that each character
occurs at most once in
$H^2(X_m)$ and those that occur are exactly those in the set
$T:=\{(b_0,b_1,b_2,b_3)\in G:\forall i:i\ne0\}\cup\{(0,0,0,0)\}$. Furthermore,
if we for $b\in\Z/m$ let $\langle b\rangle$ be the unique integer s.t.~$\langle
b\rangle\in b$ and $0\le \langle b\rangle< m$ then ${\gta}((\b
b))=1/m\sum _{i=0}^3\langle b_i\rangle$ if $(\undl b)\ne
(\undl 0)$ and 1 if not.
Finally, the action by $Gal$(\bQ/\Q) is given by $F_p((\undl
b))=(p\undl b)$ for a
prime $p\not|m$.
\begin{remark}
The proof of this in \[Ka:Sect.~6] uses transcendental methods. A purely
algebraic proof can be given by tracing the action
of $G$ through the calculations of \[SGA7:Exp.~XI].
\end{remark}
The Fermat surfaces verify the Tate conjecture (\[S-K]) so $\rk\,NS=b_2$
over a field of positive characteristic iff the average of {\gta}
over any {\gsi}-orbit equals 1. Now complex conjugations in $Gal(\bQ/\Q)$
exchanges the values 0 and 2 and fixes 1 so we see that if the subgroup
generated by {\gsi} contains complex conjugation this is always the case.
Hence if $-1\in\langle p\rangle\subseteq(\Z/m)^\times$ then
$rk~NS(X_m)=b_2(X_m)$ in characteristic p.
\begin{example}
i) $p\equiv -1 \pmod m$. Then all the orbits are of type $(1,1,\ldots,1)$
or (0,2) giving a contribution of 0 resp.~1 to ${\gsi}_0$. It is a general fact
(true whenever $rk\,NS=b_2$ and $b_2=\sum _{i+j=2}dim\;\hod Xijk$) that
$p^{2p_g}|disc\ NS$ as the morphism $H^2(X,p)\to \hod X20k$ is surjective
and vanishes on $NS\bigotimess \W$.

ii) m=5, $p\equiv 2,3 \pmod5$. Then there are four orbits of
type $(0,1,2,1)$ and
the rest are of type $(1,1,\ldots,1)$. Now the algorithm applied to $(0,1,2,1)$
gives partial sums $(-1,-1,0,0)$ and a contribution of $2$
to ${\gsi}_0$ for each copy
of this orbit and hence $disc\,NS=p^{16}$,
whereas $p_g=4$ so that we get a higher power of $p$ than is guaranteed by
$p^{2p_g}|disc\ NS$.
\end{example}
\mdef\tx{\tilde X}     \def\Tx{\tilde X}
\begin{proposition}3.1.
Suppose that $X$ is a smooth surface over \k\ and that
$p>0$. Suppose that $G$ is a finite group of order prime to
$p$ acting on $X$. Let
\tx be a minimal resolution of $X/G$. Then
$H^2(\Tx,p)=H^2(X,p)^G\perp E$, where $E$ is the \W-module
spanned by the Chern classes of the exceptional curves of $\Tx\to X/G$ and
orthogonality is wrt the cup product. Furthermore, the cup product pairing
restricted to $E$ is perfect.
\pro Let ${\gpi}\co X'\to X$ be a $G$-equivariant
blowing up of $X$ such that we have a map
${\grh}\co X'\to\Tx$ covering the quotient map $X\to X/G$. The cup product
pairing on $E$ is perfect because the cokernel of $E\to \check E$ equals \W
tensored with the sum of the local Picard groups of the singularities of $X/G$
(cf.~\[Li:14.4]) and these are killed by the order of $G$ by the
existence of a norm map. Hence we may write
$H^2(\Tx,p)=V\perp E$. Now ${\grh}_*{\grh}^*=|G|$ so ${\grh}^*$
is injective on $H^2(\tx,p)$
and the image is contained in the $G$-invariants and
is a direct factor.
Furthermore, by the projection formula, ${\grh}^* V$ is orthogonal to the
submodule of $H^2(X',p)$
spanned by the curves exceptional for {\gpi}. Therefore
${\grh}^* V\subseteq \pi^* H^2(X,p)^G$ and we are finished if we can
show that this is an equality. First, we show
this for  the $p$-torsion. Indeed, consider the slope
spectral sequence for $X'$
and \tx\ (cf.~\[Ill:II,3]). By duality
(cf.~\[Ek1]) ${\grh}_*$ is defined as a map of
spectral sequences and  we still have ${\grh}_*{\grh}^*=|G|$.
Furthermore, ${\grh}^*$ is an isomorphism on
$H^*(\tx,\W\ko \tx)\to H^*(X',\W\ko {X'})^G$ as
$H^*(X',\W\ko {X'})^G=H^*(X,\W\ko X)^G= H^*(X^G,\W\ko {X^G})
=H^*(\tx,\W\ko \tx)$ the last as the singularities are rational. Hence as
$H^0(-,\W{\gOm}_X^2)$ is torsion free as $W{\gOm}_X^2$ is we see that we
have equality on torsion groups if we have equality for the torsion of
$H^1(-,\W{\gOm}^1_-)$. The nilpotent torsion
(cf.~\[Ek1:IV,3.3.13]) of it is dual to the
nilpotent torsion of $H^2(-,\W\ko -)$
(loc.~cit.) and is hence taken care of, whereas the
semi-simple torsion (loc.~cit.:IV,3.4)
comes from $H^2(-,\Z_p(1))$ which in turn
comes from the N{\'e}ron-Severi group \[Ill:II,5.8.5] which is taken care of by
noting that
$\undl {Pic}^{\gta}(\Tx)=\undl {Pic}^{\gta}(X/G)$ as the
singularities are rational and
$\undl {Pic}^{\gta}(X/G)=\undl {Pic}^{\gta}(X)^G$ outside of the order
of $G$. Hence, as ${\grh}^* V$ is a direct factor of $\pi^* H^2(X,p)^G$,
it suffices to show that that they have the same rank. As
the rank of $V$ is the
rank of $H^2(\tx,p)$
minus the number of exceptional curves we may replace $p$ by
\l\ and then $H^2(X,\ell)^G=H^2(X/G,\ell)$
and the latter space is isomorphic to the
orthogonal complement of the exceptional curves of $\Tx\to X/G$ by the
Leray spectral sequence.
\end{proposition}

Using the proposition we get a description of the crystalline cohomology of
the minimal resolution of the quotient of $X_m$ by any subgroup of $G$.

We can also compute other invariants of Fermat surfaces and their
quotients. Consider for instance the formal Brauer group of a surface $X$
which is Mazur-Ogus (cf.~\[Ek:IV,1.1]) (in positive characteristic) or
rather $H^2(X,\W\ko X)$ the knowledge of which is equivalent to knowing the
formal Brauer group.  It follows from \[loc.~cit.:III,Thm~4.3] that
$H^2(X,W\ko X)$ is the quotient of $H^2_{cris}(X/W)\bigotimess _{W[F]}D$,
where $D$ is the Dieudonn{\'e}-ring (with power series in $V$), by the
submodule generated by $m\bigotimess 1-V(n\bigotimess 1)$ for all $m,n\in
H^2_{cris}(X/W)$ for which $Fm=pn$. Hence if $X$ is of CM-type in degree 2
we get a description of $H^2(X,W\ko X)$.
\begin{example}
i) $(0,1,0,2,1,2)$ gives a $D$-module with generators $a$ and $b$ and relations
$Fa=Vb$ and $Fb=0$. This is the Dieudonn{\'e}-module of a 2-dimensional formal
group isogenous but not isomorphic to $W_2$. For $p\equiv3\pmod7$ this
appears in the cohomology of the Fermat surface of degree 7 (the orbit of
$(1,1,1,4)\in (\Z/7)^4$).

ii) $(0,0,2,2)$  gives a $D$-module with generator a and
relation $F^2a=0$. This
is the Dieudonn{\'e}-module of $W_2$. For $p\equiv5\pmod7$ this
appears in the cohomology of the Fermat surface of degree 13 (the orbit of
$(3,3,3,4)\in (\Z/7)^4$).
\end{example}

\begin{bibliography}
\[B-Og]:\by P. Berthelot, A. Ogus \book Notes on crystalline cohomology\publ
Princeton Univ. Press \publadr Princeton\yr1978

\[C-R]: \by C. W. Curtis, I. Reiner\book Representation theory of finite
groups \yr1962\publ Interscience publishers\publadr New York

\[De]:\by P. Deligne\paper Applications de la formule des traces aux
sommes trigonom{\'e}triques \sln 569 \pages 168--232

\[D-I]: \by F. DeMeyer, E. Ingraham \paper Separable algebras over
commutative rings\sln181

\[Ek]:\by T. Ekedahl\book Diagonal complexes and F-gauge structures\publ
Hermann \publadr Paris \yr1986

\[Ek1]:\by T. Ekedahl \paper On the multiplicative properties of the
de Rham-Witt complex I\jour Arkiv f{\"o}r mate\-ma\-tik \vol22
\yr1984\pages185--239

\[Fa]:\by G. Faltings \paper p-adic Hodge theory \yr1988 \jour Journal
of the AMS \vol1\pages255--299

\[Gr]: \by M. Gros \paper Classes de Chern et classes de cycles en
cohomologie de Hodge-Witt logarithmique\jour M{\'e}m. de la Soc. Math. de
France\vol21\pages1--87\yr1985

\[Gro]:\by A. Grothendieck \paper Le groupe de Brauer III \inbook Dix
expos{\'e}s sur la cohomologie des sch{\'e}mas \publ North-Holland \publadr
Amsterdam \yr1968 \pages88--188

\[Ill]:\by L. Illusie \paper Complexe de de Rham-Witt et cohomologie
cristalline \jour Ann. scient. {\'E}c. Norm. Sup \vol 12 \yr1979
\pages501--661

\[Jo]:\by J. de Jong\spec Article to appear

\[Ka]:\by N. M. Katz\paper On the intersection matrix of a hypersurface
\jour Ann. scient. {\'E}c. Norm. Sup\vol2\yr1969\pages 583--589

\[K-L]:\by N. M. Katz, S. Lang \paper Finiteness theorems in geometric
class field theory \jour L'Ens. Math \vol27 \yr1981\pages286--314

\[K-M]: \by N. M. Katz, W. Messing \paper Some consequences of the Riemann
hypothesis for varieties over finite fields \jour Invent. Math. \yr1974
\vol23\pages73--77

\[Li]:\by J. Lipman\paper Rational singularities with applications to
algebraic surfaces and unique factorization \jour Publ. IHES
\vol36\yr1969\pages195--280

\[Re]:\by I. Reiner\book Maximal orders\publ Academic Press\publadr London
\yr 1975

\[Se]:\by J.-P. Serre \book Abelian \l-adic representations and elliptic
curves \publ W. A. Benjamin\publadr New York\yr1968

\[Se1]:\by J.-P. Serre \paper Groupes alg{\'e}briques associ{\'e}s aux modules
de Hodge-Tate \jour Ast{\'e}rix 65 \yr 1979\pages155--188

\[SGA7]: \by P. Deligne, N. M. Katz \paper Groupes de monodromie en
g{\'e}om{\'e}trie alg{\'e}brique \sln340

\[S-K]:\jour Toh{\^o}ku Math. Jour. \by T. Shioda, T. Katsura \yr 1979
\vol31\paper On Fermat varieties \pages 97-115

\[S-I]: \by T. Shioda, H. Inose \paper On singular K3-surfaces \inbook
Complex analysis \& algebraic geometry \publ Cambridge Univ. Press
\publadr Cambridge \yr1977\pages 117--136

\end{bibliography}
\end{section}
\end{document}